\documentclass[11pt]{amsart}
\usepackage{amssymb}
\usepackage{amsfonts}
\usepackage{amscd}
\usepackage{amsmath}
\usepackage{mathrsfs}
\usepackage{curves}
\usepackage{epsfig}
\usepackage[pass]{geometry}
\usepackage{graphicx}
\usepackage{verbatim}
\usepackage{latexsym}
\usepackage{eucal}
 
\newtheorem{theorem}{Theorem}[section]
\newtheorem{lemma}[theorem]{Lemma}
\newtheorem{prop}[theorem]{Proposition}

\newtheorem{remark}[theorem]{Remark}

\def \mca {{\mathcal A}}

\def \mcd {{\mathcal D}}
\def \mce {{\mathcal E}}
\def \mcf {{\mathcal F}}
\def \mcg {{\mathcal G}}

\def \mcl {{\mathcal L}}

\def \mcp {{\mathcal P}}

\def \mcr {{\mathcal R}}

\def \mcu {{\mathcal U}}
\def \mcv {{\mathscr V}}
\def \mcw {{\mathcal W}}

\def \msc {{\mathscr C}}

\def \msm {M}
\def \msx {{\mathcal X}}
\def \msv {{\mathscr V}}

\def \mbc {{\mathbb C}}

\def \mbn {{\mathbb N}}
\def \mbr {{\mathbb R}}
\def \mbs {{\mathbb S}}

\def \id {\operatorname{Id}}

\def \im {\operatorname{Im}}

\def \supp {\text{supp }}

\def \beqq {\begin{equation}}
\def \eeqq {\end{equation}}

\def \bpf {\begin{proof}}
\def \epf {\end{proof}}
\def \beq {\begin{equation*}}
\def \eeq {\end{equation*}}

\def \eps {\epsilon}   
\def \la {\lambda}   
\def \La {\Lambda}    
   
\def \lap {\Delta}
\def \p {\partial}
\def \ha {\frac{1}{2}}

\def \tilde {\widetilde}

\def \op {\text{Op}}

\begin{document}
\title[]{The Anisotropic Calder\'on Problem for High Fixed Frequency}
\author{Gunther Uhlmann}
\address{Gunther Uhlmann
\newline
\indent Department of Mathematics, University of Washington 
\newline
\indent and Institute for Advanced Study, the Hong Kong University of Science and Technology}
\email{gunther@math.washington.edu}
\author{Yiran Wang}
\address{Yiran Wang
\newline
\indent Department of Mathematics, Emory University}
\email{yiran.wang@emory.edu}
\begin{abstract}
We consider Schr\"odinger operators at a fixed high frequency on simply connected compact Riemannian manifolds with non-positive sectional curvatures and smooth strictly convex boundaries. We prove that the Dirichlet-to-Neumann map uniquely determines the potential. 
\end{abstract}
\date{\today}
\maketitle

 \section{Introduction}
Let $(\msm, g)$ be a simply connected compact Riemannian manifold of dimension $n = 3$ with smooth strictly convex boundary $\p \msm$.   Let $\lap_g$ be the positive Laplace-Beltrami operator on $(\msm, g)$. For $V\in C_0^\infty(\msm)$ and $\la >0$, we consider the Dirichlet problem 
\beqq\label{eq-dtn0}
\begin{gathered}
(\lap_g + V - \la^2) u  = 0  \text{ in } \msm, \\
u = f \text{ at } \p \msm
\end{gathered}
\eeqq
Suppose that $\la^2$ is not an eigenvalue of $\lap_g + V.$ Let $u$ be the unique $C^\infty(\msm)$ solution for $f\in C^\infty(\p \msm)$.  
We consider the Dirichlet-to-Neumann map $\La$ defined by 
\beqq\label{eq-dtn1}
\La f \doteq \p_\nu u = |g|^\ha  \sum_{i, j = 1}^3 g^{ij} (\p_i u) \nu_j
\eeqq
 where $\nu$ is the unit outward normal to $\p \msm$. The problem we study is the determination of $V$ from $\La$ for a large but fixed $\la$. Our main result is 
\begin{theorem}\label{thm-main}
Let $(\msm, g)$ be a simply connected compact Riemannian manifold of dimension three with smooth boundary $\p \msm$ which is strictly convex. Assume that the sectional curvatures are non-positive.  Let $V, \tilde V \in C_0^\infty(\msm)$ be two potentials which are supported away from $\p M$, and $\La, \tilde \La$ be the corresponding Dirichlet-to-Neumann map for \eqref{eq-dtn0}. Suppose that $\la$ is sufficiently large and $\la^2$ is not an eigenvalue for \eqref{eq-dtn0} for both $V$ and $\tilde V$. Then $\La = \tilde \La$ implies  $V = \tilde V.$ 
\end{theorem}

In addition to the uniqueness result, we remark that for a class of potentials, it is possible to obtain a H\"older type stability estimate from our proof. Also, our analysis leads to an approximate reconstruction method for high frequencies by inverting a weighted geodesic ray transform on $(M, g)$. These will be discussed in the end of Section \ref{sec-proof}.

When $\la=0$ this problem is related to the anisotropic Calder\'on problem. In fact the Calder\'on problem in the isotropic case can be reduced to study an inverse boundary problem for the Schr\"odinger equation at zero energy. It is well known that in dimension three or larger the anisotropic Calder\'on problem can be reduced to study the Dirichlet-to-Neumann map for the Laplace-Beltrami operator of a Riemannian metric uniquely determined by the conductivity,  see  \cite{LU}. The case we are considering here corresponds to determining a metric in the same conformal class. For the Euclidean metric the problem we are considering here at any fixed energy was solved in dimension three or larger by Sylvester and Uhlmann \cite{SyUh} and in two dimensions by Bughkeim \cite{Bug}. For the anisotropic case on conformally transversally anisotropic manifolds this problem was solved under some conditions on the transversal manifold, see  \cite{FKSU},   \cite{FKLS}. In two dimensions using Bughkeim's method the problem was solved on Riemann surfaces in \cite{GuTz}. The problem of recovering the metric itself has only been solved in dimensions three or larger in the real-analytic category see  \cite{LaUh}  and \cite{LTU}. We plan to consider the problem of recovering the Riemannian metric at large fixed frequency in a subsequent article. For further references and a review of Calder\'on's problem, see  \cite{Uh}.

We will outline the proof of Theorem \ref{thm-main} in Section \ref{sec-stra}, but we want to point out that some assumptions in the statement of the theorem are made to simplify  the analysis in order to give a clear demonstration of our method. They can potentially be removed. First of all, the parametrix construction we use in the proof is simpler in dimension three with the non-positive curvature assumptions, but it can be done for higher dimensions (and even without the curvature assumptions). Second, we will use stability estimates for geodesic ray transforms on $(\msm, g)$ which is simple to state with the non-positive sectional curvature assumptions. It is known that such estimates hold under weaker curvature conditions, see  \cite{Sha}. In fact, the stability estimates are already known for manifolds satisfying the foliation condition, see  \cite{UV}. Finally, we assumed that $V, \tilde V$ are supported away from $\p\msm$ to save us from some technical discussions related to the singularities of the distance function.   
 
 \section{The strategy}\label{sec-stra}
 As we concern high frequencies, we can assume $\la>> 0$ and take $h = 1/\la$ as a semiclassical parameter. We consider the semiclassical problem 
\beqq\label{eq-semi}
\begin{gathered}
h^2(\lap_g + V)  u - \sigma^2 u  = 0  \text{ in } \msm, \\
 u = f \text{ at } \p \msm
\end{gathered}
\eeqq 
for $h\in (0, 1]$. Here, we allow $\sigma$ in a compact set of $\mbc\backslash 0$ but $\sigma$ will be fixed in later arguments. We will use Green's representation formula which connects $\La$ to the fundamental solution of the  equation in \eqref{eq-semi}. 

Since $(\msm, g)$ has non-positive sectional curvature and $M$ is simply connected, without loss of generality we  take $\msm$ as a subset of $\mbr^3.$ Let $\mcr(\sigma, h)$ be the fundamental solution such that
\beqq\label{eq-funda}
(h^2(\lap_g + V) - \sigma^2)\mcr(\sigma, h) = \id \text{ in } \mbr^3
\eeqq
Here, we can extend $V$ to be zero outside $M$ so $V$ is a $C_0^\infty$ potential on $\mbr^3$.  The construction of the fundamental solution is known, see  e.g.\ \cite{Shi}. In particular, for fixed $h, \sigma$ and modulo a smoothing operator, $\mcr(\sigma, h)$ is a pseudo-differential operator of order $-2$ on $\mbr^3$.   

Next, we follow the argument in  \cite[page 80]{Tre}. Let $u$ be the solution of \eqref{eq-semi}. 
Let $\chi_\msm$ be the characteristic function for $\msm$ in $\mbr^3$. Then we have 
\beq
  \lap_g (\chi_\msm u) =   \chi_\msm \lap_g u  + 2   \sum_{i, j = 1}^n g^{ij} \frac{\p u}{\p z^i} \frac{\p \chi_\msm}{\p z^j} + u \lap_g \chi_\msm
\eeq
Thus, 
\beq
  h^2(\lap_g + V) (\chi_\msm u)  - \sigma^2 \chi_\msm u  =    2  h^2 \sum_{i, j = 1}^n g^{ij} \frac{\p u}{\p z^i} \frac{\p \chi_\msm}{\p z^j} +  u  h^2(\lap_g \chi_\msm).
\eeq
Using the fundamental solution $\mcr(\sigma, h)$, we obtain for $z\in \msm$ that 
\beqq\label{eq-dist1}
\chi_\msm u =    2 h^2 \mcr(\sigma, h) \ast (\sum_{i, j = 1}^n g^{ij} \frac{\p u}{\p z^i} \frac{\p \chi_\msm}{\p z^j}) +  \mcr(\sigma, h) \ast( u (h^2 \lap_g \chi_{\msm})).
\eeqq
Following the calculation in  \cite[page 80]{Tre}, we get that 
\beqq\label{eq-ide0}
\begin{gathered}
u(z) 
=   -   h^2\int_{\p \msm} \mcr(z, z', \sigma, h) \frac{\p u}{\p \nu}(z') d\sigma(z') 
+  h^2\int_{\p \msm} u(z') \frac{\p }{\p \nu}  \mcr(z, z', \sigma, h) d\sigma(z') 
\end{gathered}
\eeqq
In this work, we will use the same notation $\mca$ for both the operator $\mca: \mce'(\msm)\rightarrow \mcd'(\msm)$ and its Schwartz kernel. 
Formula \eqref{eq-ide0} holds in the interior of $\msm$. But we will use this formula for $z\in \p \msm$ which  will be justified in Section 4. 

Now suppose we have  two potentials $V, \tilde V$ on $\msm$. Let $\La, \tilde \La$ be the corresponding DtN maps. Let  $\mcr, \tilde \mcr$ be the corresponding fundamental solutions. For $(z, z')\in \p \msm\times \p\msm$ away from $z = z'$,  we use \eqref{eq-ide0} to get
\beqq\label{eq-ide}
\begin{gathered}
 \p_\nu \tilde \mcr(z, z', \sigma, h) - \p_\nu \mcr(z, z',\sigma, h) \\
 = \int_{\p \Omega} \tilde \mcr(z, z'', \sigma, h) \tilde\La(z'', z') dz''   - \int_{\p \Omega} \mcr(z, z'', \sigma, h)  \La(z'', z') dz''  
\end{gathered}
\eeqq 
Our proof is based on the investigation of this formula. We list the key steps and describe the structure of the paper. 
 
{\bf Step 1:} We construct an approximation of the kernel of $\mcr(\sigma, h)$ with an explicit leading order term as $h\rightarrow 0$. This is done in Section \ref{sec-para} and Section \ref{sec-res}. Roughly speaking, $\mcr(\sigma, h) = \mcg(\sigma, h) + O(h^{-1})$ and we have 
\beq
\tilde \mcr(\sigma, h) -  \mcr(\sigma, h)  = \tilde \mcg(\sigma, h) - \mcg(\sigma, h) + \mcf_{para}
\eeq
where $\mcf_{para}$ denotes the error term.  Here, $\mcg$ is of the form 
\beqq\label{eq-gform}
h^{-2}e^{-i\frac\sigma h r(z, z')} r(z, z')^{-1}\mca(z, z'), \quad z, z'\in \msm
\eeqq
where $\mca$ is smooth and $r$ denotes the distance function on $\msm$. Because there is no conjugate points on $\msm$, the distance function is smooth on $\msm\times \msm$ away from the diagonal $\{(z, z')\in \msm\times \msm: z = z'\}$.

{\bf Step 2:} We show that the leading order term in $\tilde \mcg(\sigma, h) - \mcg(\sigma, h)$ as $h\rightarrow 0$ is a weighted geodesic ray transform of the difference $W = \tilde V - V$, denoted by $\msx^w W$ so
\beq
\tilde \mcg(\sigma, h) - \mcg(\sigma, h)  = \msx^w W + \mcf_{lead, 1}
\eeq
where $\mcf_{lead, 1}$ denotes the error term. This is obtained from some perturbation argument in Section \ref{sec-pert} and a stationary phase argument by considering the oscillatory behavior in \eqref{eq-gform} done in Section  \ref{sec-geo}. 
The error term  comes from the application of stationary phase argument. It is smaller in terms of $h$, but involves higher order derivatives of $W$. 

{\bf Step 3:} For $\msx^w W$, we use the stability estimate in  \cite{Sha}. After some careful estimates of the error terms which is done in Section \ref{sec-est} and semiclassical analysis of the DtN map in Section \ref{sec-dtn} and \ref{sec-heat}, we derive  that  
\beqq\label{eq-eq1}
\|W\|_{L^2(\msm)} \leq C  h |W|_{C^2(\msm)}
\eeqq
where $C$ is a generic constant.  Here, for $u\in C^m(\msm),$ we denote the seminorm $|u|_{C^m(\msm)} = \sup_{x\in \msm} \sum_{|\alpha|\leq m} |\p^\alpha u|$. For given $V, \tilde V$, there exists a constant $C_1$ such that $|W|_{C^2(\msm)}\leq C_1 \|W\|_{L^2(\msm)}$. 
We thus obtain that $W = 0$ for $h$ sufficiently small (depending on $V, \tilde V$). This is done in Section \ref{sec-proof}.

\section{The semiclassical parametrix}\label{sec-para}
Consider the semiclassical operator 
\beq
\mcp \doteq h^2 (\lap_g + V) - \sigma^2
\eeq
 on $\mbr^3$ where $V\in C_0^\infty(\mbr^3)$. Denote the resolvent by 
 \beq
 \mcr(\sigma, h) = (h^2(\lap_g  + V) - \sigma^2)^{-1}.
 \eeq  
We follow the discussion in  \cite[Page 24]{MSV} to construct an approximation of $\mcr(\sigma, h)$.

Let $r$ be the distance function on $(\msm, g)$.  In polar coordinates $(r, \theta)$ based at a point $z'\in \msm$, the metric 
\[
g = dr^2 + H(r, \theta, d\theta)
\]
where $H$ is a smooth $1$-parameter family of metrics on $\mbs^2.$ This is Gauss lemma, see \cite[page 91]{GHL}. The Laplacian in this coordinate reads
\[
\lap_g  = -\p_r^2 - A\p_r + \lap_H, \quad A = |g|^{-\ha}\p_r (|g|^\ha)
\]
where $|g|^\ha$ is the volume element  and $\lap_H$ is the positive Laplacian on $\mbs^2$ with respect to $H(r, \theta, d\theta)$. 
 As $r\rightarrow 0$, we have
\[
|g|^{\ha}(r, \theta) = r^2(1 + r^2 g_1(r, \theta)),
\]
where $g_1$ is smooth at $r = 0$, see \cite[page 144]{GHL}. Thus 
\[
\lap_g = -\p_r^2 - (2/r+ r A(r, \theta))\p_r + \lap_H
\]
in which $A(r, \theta)$ is smooth up to $r = 0.$ 

Now we look for an approximation of the resolvent whose Schwartz kernel is of the form 
\[ 
\mcg(\sigma, h, z, z') = e^{-i   \frac{\sigma}{h}r(z, z')} (h^{-2} \mcu_0(z, z') + h^{-1}\mcu_1(z, z')) 
\]
We formally compute that  
\beqq\label{eq-gasymp}
\begin{gathered}
(h^2 ( \lap_g + V)  - \sigma^2) \mcg(\sigma, h) = e^{-i \frac{\sigma}{h}r} (2i  \sigma h^{-1} |g|^{-1/4}\p_r (|g|^{\frac 14}\mcu_0)  \\
+ 2i \sigma |g|^{-1/4}\p_r (|g|^{1/4}\mcu_1) + (\lap_g + V) \mcu_0  + h(\lap_g + V) \mcu_1
\end{gathered}
\eeqq
 First, the $h^{-1}$ terms have to vanish and this gives
\[
 |g|^{-1/4}\p_r (|g|^{\frac 14}\mcu_0) = 0, \quad r>0
\]
Notice that $|g|^\ha$ is the density factor on $\msm\times \msm$ so $|g|^{1/4}$ is indeed the half-density factor and the above equation is the Lie derivative of the principal symbols, see \cite{MSV}. We solve that 
\beq
\begin{gathered}
\mcu_0 = \frac{1}{4\pi} |g|^{-1/4}  
\end{gathered}
\eeq
Near $r = 0$, we get 
\beq
\mcu_0 =  \frac{1}{4\pi}(r^{-1} + rB)
\eeq
where $B$ is smooth in $r, \theta$ up to $r = 0$. This implies that 
\[
\lap_g \mcu_0 = \delta(z, z') + \frac{1}{4\pi} A r^{-1} + \frac{1}{4\pi} \lap_g(rB). 
\]
So we removed the singularities at the diagonal to the leading order. This only happens in dimension three. 

Next, the $h^0$ terms in \eqref{eq-gasymp} have to vanish and we obtain
\beq
\begin{gathered}
2i \sigma |g|^{-1/4}\p_r (|g|^{1/4}\mcu_1) +  (\lap_g + V) \mcu_0 = 0, \quad r > 0\\
\mcu_1 = 0 \text{ at } r = 0
\end{gathered}
\eeq
So we get 
\[
\mcu_1 = -\frac{1}{2 i \sigma} |g|^{-1/4} \int_0^r |g|^{\frac 14} (\lap_g + V) \mcu_0 ds.
\]
In particular, $|g|^{1/4}$ is smooth up to $r = 0$ and vanishes at $r = 0$. Thus the integrand in $\mcu_1$ is smooth at $r = 0.$ This implies that $\mcu_1$ is smooth up to $r = 0. $ 

Consider the remainder term. We have
\beq
\begin{gathered}
(h^2 (\lap_g + V) - \sigma^2) \mcg(\sigma, h)  = \delta(z, z') + h  e^{-i \frac{\sigma}{h}r} (\lap_g + V) \mcu_1\\
= \id + h \mce(\sigma, h)
\end{gathered}
\eeq
We observe that the remainder term $\mce$ has a $1/r$ type singularity at $r = 0.$ Thus, in this approach, we removed the semiclassical error to order $h$ but the in the classical sense, we only removed the leading order singularity at $r = 0$. 

 We are done with the construction and we summarize the result.  
\begin{prop}\label{prop-para}
For $h\in (0, 1], \sigma \in \mbc\backslash 0$, there exists an operator $\mcg(\sigma, h)$ and $\mce(\sigma, h)$ such that  
\beq
(h^2 (\lap_g + V) - \sigma^2) \mcg(\sigma, h) = \id + h \mce(\sigma, h)
\eeq
where the Schwartz kernel 
\[
\mcg(\sigma, h, z, z') = e^{-i\frac{\sigma}{h}r(z, z')} (h^{-2} \mcu_0(z, z') + h^{-1} \mcu_1(z, z'))
\]
with 
\beq
\mcu_0 =  \frac{1}{4\pi} |g|^{-1/4}, \quad \mcu_1 =  -\frac{1}{2\pi i \sigma} |g|^{-1/4} \int_0^r |g|^{\frac 14} (\lap_g  + V ) \mcu_{0} ds
\eeq
and 
\beq
\mce(\sigma, h, z, z') =  e^{-i \frac{\sigma}{h}r(z, z')} (\lap_g + V) \mcu_1(z, z').
\eeq 
 \end{prop}
 
 \begin{remark}
 For higher dimensions $n\geq 4$, a similar but more involved construction was also done in \cite{MSV} for small metric perturbations of hyperbolic spaces and it was further developed in \cite{SaWa} for non-trapping asymptotically hyperbolic manifolds.  For our purpose, the asymptotic behavior near the infinity does not matter. The construction in \cite{MSV} and \cite{SaWa} away from the infinity does not rely on the hyperbolic structure and can be modified to obtain a parametrix for simply connected manifolds. \end{remark}

\section{The resolvent kernel and estimates}\label{sec-res}
Hereafter, we fix $\sigma \in \mbc\backslash 0$. We will drop $\sigma$ in the notations and write e.g.\ $\mcr(h) = \mcr(\sigma, h)$. We consider the resolvent $\mcr(h)$ and  use the parametrix to find the resolvent kernel. The construction and estimate for the Green's function are known for elliptic problems in general, see for example \cite[Chapter VI, Section 4]{Shi}. 
Here, the point is  the dependency on $h$. 
 
We start from 
\beqq\label{eq-res1}
(h^2 (\lap_g + V) - \sigma^2) \mcg(h) = \id + h \mce(h)
\eeqq
Using the resolvent $\mcr(h)$ in \eqref{eq-res1}, we get 
\beq\label{eq-reside0}
\mcg(h) = \mcr(h)(\id +  h \mce(h))  
\eeq
We first have 
\begin{lemma}
Let $\sigma \in \mbc \backslash 0$ and $\im \sigma \leq 0.$ 
There is $h_0> 0$ such that for $0 < h < h_0$, the operator $\id + h \mce(h)$ is invertible on $L^2(\msm)$.  
\end{lemma}
\bpf
We recall that 
 \beq
 r\lap_g  =  -r \p_r^2 - (2+ r^2 A(r, \theta))\p_r + r\lap_H
 \eeq
 is a second order differential operator  with smooth coefficients up to $r = 0.$ From the formula of $\mce(z, z', h)$ in Proposition \ref{prop-para}, we see that $r(z, z')\mce(z, z',  h)$ is smooth up to $r = 0.$  
 We use Schur's lemma to estimate the $L^2$ norm of $\mce(h)$. First, we have
\beq
 \int_{\msm} |\mce(h, z, z')| dz'   \leq C \frac{1}{|\sigma|}  \int_\msm \frac{1}{r(z, z')}dz' \leq C \frac{1}{|\sigma|} 
\eeq
where the constant $C$ depends on $\msm$ but not on $h, \sigma$. Here, we used $\im \sigma \leq 0$. Using Schur's lemma,  we get 
\beq
\|h\mce(h)\|_{L^2(\msm)\rightarrow L^2(\msm)} \leq C  h|\sigma|^{-1}
\eeq
For $h_0 > 0$ such that for $h < h_0|\sigma|^{-1}$, we have $C  h|\sigma|^{-1} < 1$ and thus  $\id + h \mce(h)$ is invertible on $L^2(\msm)$.  
\epf

The proof implies that for $h <h_0,$ the inverse can be written as a Neumman series. We write
\beqq\label{eq-reside1}
\begin{gathered}
\mcr(h) = \mcg(h)(\id + h \mce(h))^{-1} = \mcg(h) - \mcg(h) h \mce(h) + \cdots\\
 = \mcg(h) + \mcf(h), 
 \text{ where } \mcf(h) = \sum_{j = 1}^\infty \mcg(h)(-h \mce(h))^j
 \end{gathered}
\eeqq
We find out the Schwartz kernel of $\mcf(h).$  
We start with some estimates of the terms in the parametrix.  
 \begin{lemma} 
There exists $C>0$ depending on $g, |V|_{C^0(\msm)}$ such that for fixed $\sigma\in \mbc\backslash 0$ and $\im \sigma\leq 0,$  we have
\begin{enumerate}
\item $|r(z, z')\mcg(h,  z, z')|_{C^0(\msm\times \msm)}   \leq  Ch^{-2}  $
 \item $|r^2(z, z')\p_\nu \mcg(h, z, z')|_{C^0(\msm\times \msm)}   \leq C  h^{-3}  $
 \item  $|r(z, z')\mce(h,  z, z')|_{C^0(\msm\times \msm)} \leq C  $
\end{enumerate}
\end{lemma}
\bpf 
 This is straightforward. For (1), the expression of $\mcg$ is 
 \[
\mcg(h, z, z') = e^{-i\frac{\sigma}{h}r(z, z')} (h^{-2} \mcu_0(z, z') + h^{-1} \mcu_1(z, z'))
\] 
From the expression of $\mcu_0, \mcu_1$, we know that $r \mcu_0, r \mcu_1$ are both smooth up to $r = 0.$ 
Therefore, we get that 
\beq
\begin{gathered}
 |r(z, z')\mcg(z, z', h)|  \leq Ch^{-2} + Ch^{-1} |\sigma|^{-1}
\end{gathered}
\eeq
for some constant only depends on $g, |V|_{C^0(\msm)}$.  For (2), we have 
\beq
\begin{gathered}
  \p_\nu \mcg(z, z', h) = -i \frac{\sigma}{h} \p_\nu r(z, z') e^{-i\frac{\sigma r(z, z')}{h}} (h^{-2} \mcu_0 + h^{-1} \mcu_1) 
+ e^{-i\frac{\sigma r(z, z')}{h}} \p_\nu (h^{-2} \mcu_0 + h^{-1} \mcu_1) 
\end{gathered}
\eeq
Therefore, 
\beq
\begin{gathered}
 |r^2(z, z')\p_\nu \mcg(z, z', h)|  \leq C |\sigma| h^{-3} + Ch^{-2}  + Ch^{-1}|\sigma|^{-1}
 \end{gathered}
\eeq
The proof of (3)  follows from the same argument.   
\epf

\begin{lemma}
Let $\sigma \in \mbc\backslash 0$ and $\im\sigma \leq 0$. There exists $h_0>0$, $C>0$ depending on $g, |V|_{C^0(\msm)}$ and $\sigma$ such that  for $h < h_0$, we have
\begin{enumerate}
\item $\displaystyle{|r(z, z')\mcf(h, z, z')|_{C^0(\msm\times\msm)} \leq C h^{-1}}$
\item $\displaystyle{|r^2(z, z')\p_\nu \mcf(h, z, z')|_{C^0(\msm\times\msm)} \leq C  h^{-2}}$
\end{enumerate} 
\end{lemma}
\bpf
We  estimate the Schwartz kernel of $h^j\mcg(h)\mce^j(h), j\geq 1$. 
The kernel can be written as  
\beq
\begin{gathered}
h^j\mcg(h)\mce^j(h)(z_0, z_j) = h^j\int_{\msm\times \msm\times \cdots \msm} \mcg(h, z_0, z_1) \mce(h, z_1, z_2) \mce(h, z_2, z_3)\\
\cdots \mce(h, z_{j-1}, z_j) dz_2 \cdots dz_{j-1} 
\end{gathered}
\eeq
where $z_0, z_j\in \msm.$ Since $r(z, z')\mce(h, z, z')$ are smooth up to $r = 0$,  we can estimate
\beq
\begin{gathered}
|h^j \mcg(h)\mce^j(h)(z_0, z_j)| \leq   C^j h^j |\sigma|^{-j} \int_{\msm\times \msm\times \cdots \msm} |\mcg(h, z_0, z_1)| \frac{1}{r(z_1, z_2)} \frac{1}{r(z_2, z_3)} \\
\cdots \frac{1}{r(z_{j-1}, z_j)} dz_2 \cdots dz_{j-1} 
\end{gathered}
\eeq
The kernel is integrable. Moreover, we observe that 
\beq
\int_{M} \frac{1}{r(x, y)} \frac{1}{r(y, z)} dy \leq \frac{1}{r(x, z)} C \int_{M} \frac{1}{r(x, y)} \frac{1}{r(y, z)} dy \leq C \frac{1}{r(x, z)}
\eeq
where $C$ only depends on $(M, g)$. So we get 
\beq
\begin{gathered}
|h^j r(z, z') \mcg(h)\mce^j(h)(z, z')| \leq   (C_0h/|\sigma|)^j Ch^{-2}
\end{gathered}
\eeq
for some $C_0 > 0$. For $h/|\sigma| < 1$, we can sum the terms from $j = 1$ and get 
\beq
\begin{gathered}
|r(z, z')\mcf(h, z, z')|_{C^0(\msm\times \msm)} \leq \frac{C_0 h/|\sigma|}{1 - C_0 h /|\sigma|} Ch^{-2} 
\leq Ch^{-1}|\sigma|^{-1}  
\end{gathered}
\eeq
where we used $h/|\sigma| < 1/2$.  The derivative estimate is similar. 
\epf
 
 \begin{remark}
Now we justify the formula \eqref{eq-ide0} for $z\in \p \msm$.  First of all, the first term on the right of \eqref{eq-ide0} is continuous at $\p \msm$ because the kernel $\mcr(z, z', \sigma, h)$ is locally integrable.  For the last term in \eqref{eq-ide0}, we notice that the kernel $\p_\nu \mcr(z, z', \sigma, h)$ is not locally integrable but the kernel is only singular at $z = z'$. We will stay away from the diagonal as follows. For any $p \in \p \msm$, let $\chi_{p}$ be a compactly supported smooth cut off function supported near $p$. For any $f\in C^\infty(\p \msm),$ we consider Dirichlet data $f\chi_p$. Then the last term is smooth when we consider $z$ outside of support of $\chi_{p}$. This means that if we stay away from the diagonal, all the terms in \eqref{eq-ide0} can be extended continuously  to $\p \msm$.  This is how we get \eqref{eq-ide}. 
\end{remark}

\section{The perturbation argument}\label{sec-pert}
Let $V, \tilde V$ be two potentials on $(\msm, g)$. 
 Let $\mcr(h), \tilde \mcr(h)$ be the corresponding resolvent of $h^2(\lap_g + V) -\sigma^2, h^2(\lap_g + \tilde V) -\sigma^2$. We are interested in the difference $\tilde \mcr(h) - \mcr(h)$. Let $W = \tilde V - V$. From the resolvent formula, we get
 \beq
 \tilde \mcr(h) - \mcr(h) = \tilde \mcr(h)h^2W \mcr(h)
 \eeq
 so that 
 \beqq\label{eq-resdiff}
 \begin{gathered}
 \tilde \mcr(h) = \mcr(h)(\id + h^2W\mcr(h))^{-1}  
 \end{gathered}
 \eeqq
 Here, because $W$ is compactly supported in $\msm$, the invertibility of $\id + h^2W\mcr(h)$ on $L^2(\msm)$ follows from the analytic Fredholm theory, see page 19-20 of \cite{Mel}. If we apply $h^2 (\lap_g + V) - \sigma^2$ to \eqref{eq-resdiff}, we get 
 \beq
 (h^2 (\lap_g + V) - \sigma^2)\tilde \mcr(h) = (\id + h^2W\mcr(h))^{-1}  
 \eeq
 Then it follows from the structure of $\tilde \mcr$ in Section \ref{sec-res} that $(\id + h^2W\mcr(h))^{-1}$ has an integrable kernel and the $L^1$ norm of the kernel is bounded  as $h\rightarrow 0.$
 
Now we describe the approximation of \eqref{eq-ide} that we use later.   
From \eqref{eq-resdiff}, we write
\beq
\begin{gathered}
\tilde \mcr(h) -  \mcr(h) = - \mcr(h) h^2 W \mcr(h) + \mcf_{res}, 
\end{gathered}
\eeq
where 
  \beqq\label{eq-fres}
 \begin{gathered}
\mcf_{res} 
 = \mcr(h) h^2 W \mcr(h) h^2 W \mcr(h)(\id + h^2 W\mcr(h))^{-1}.
 \end{gathered}
 \eeqq
Here,  $\mcf_{res}$ accounts for the error in the potential perturbation. We use the parametrix $\mcr(h) = \mcg_{0} + \mcg_{1} + \mcf$ to get 
 \beqq\label{eq-fpara}
 \begin{gathered}
 -h^2\mcr(h)W \mcr(h)  = \mcf_{lead} + \mcf_{para}, \\
 \mcf_{lead} = -\mcg_{0} h^2 W \mcg_{0}, \\
 \mcf_{para} = -\mcg_{0} h^2 W (\mcg_{1} + \mcf) - (\mcg_{1} + \mcf) h^2 W \mcg_{0} - (\mcg_{1} + \mcf) h^2 W (\mcg_{1} + \mcf)
 \end{gathered}
 \eeqq
 The term $\mcf_{lead}$ is what we use to get a weighted geodesic ray transform of $W$. The term $\mcf_{para}$ accounts for the error in parametrix construction.  
To summarize, we get
 \beqq\label{eq-approx}
 \begin{gathered}
 \mcr(z, z',  h) - \tilde \mcr(z, z', h)  
     =  \mcf_{lead} + \mcf_{para} + \mcf_{res} 
     \end{gathered}
 \eeqq
 In the next section, we analyze $\mcf_{lead}$. Then we estimate $\mcf_{para}, \mcf_{res}$.

 \section{Geodesic ray transform}\label{sec-geo}
By the assumption that $(\msm, g)$ has non-positive sectional curvatures, we know that for every $z, z'\in \p \msm, z\neq z'$, there is a unique distance minimizing geodesic $\gamma_{z, z'}$ between them. Let $r(z, z')$ be the distance between $z, z'\in \msm$. Let $\gamma_{z, z'}(s): [0, r(z, z')] \rightarrow \msm$ be the unit speed geodesic from $z$ to $z'$. It satisfies the geodesic equation 
\beq
\begin{gathered}
\nabla_{\dot \gamma_{z, z'}(s)} \dot \gamma_{z, z'}(s) = 0\\
\gamma_{z, z'}(0) = z, \quad \dot \gamma_{z, z'}(0) =\p_z r(z, z')
\end{gathered}
\eeq
For $z, z'\in \p \msm,$ we consider the geodesic ray transform on scalar functions 
\beqq\label{eq-xray0}
\msx f(z, z') = \int_0^{r(z, z')} f(\gamma_{z, z'}(s)) ds
\eeqq
Note that we parametrize the geodesics using $z, z' \in \p \msm.$ Usually, the geodesic transform is parametrized by using inward pointing unit tangent bundle at $\p \msm$
\beq
\Omega_-\msm = \{(z, \xi)\in T\msm | z\in \p \msm, -\langle \xi, \nu\rangle \geq 0, |\xi|_g = 1\}
\eeq
 see for example \cite{Sha}. For $(\msm, g)$ that we consider, there is no conjugate points. There is a diffeomorphic between $(z, z')\in\p\msm\times\p\msm$ away from the diagonal $z = z'$ and $(z, \xi)\in \Omega_-\msm$ away from $\xi = 0$ via $\xi = \p_zr(z, z')$. For our purpose, the function $f$ in \eqref{eq-xray0} is supported away from $\p M$ thus it suffices to consider geodesics corresponding to $(z, z')\in \p\msm\times\p \msm$ away from the diagonal. For this reason, we can use $\msc = \p \msm\times \p \msm$ with a measure which away from the diagonal is the one induced from $\Omega_-\msm$. We use $\msc$ as the set for parametrizing geodesics.
 
 We   recall the following stability estimates of the geodesic ray transform which is  a restatement of Theorem 4.3.3 of \cite{Sha} for our geometric setting. 
\begin{theorem}\label{thm-xraystab}
Let $(\msm, g)$ be a simply connected compact Riemannian manifold with strictly convex boundary $\p \msm$ and non-positive sectional curvatures. Suppose $f\in H^1(\msm)$ and $f$ is supported away from $\p \msm.$ Then $f$ is uniquely determined by the ray transform $\msx f$ and the conditional stability estimate holds
\beq
\|f\|_{L^2(\msm)} \leq C \|\msx f\|_{H^1(\msc)}
\eeq 
where $C$ is a constant independent of $f$.
\end{theorem}

Later, we shall consider $W\in C^\infty(\msm)$ and $\mcw$ a smooth weight function  on $\msm$ with the property $C_1\leq |\mcw|_{C^0(\msm)}\leq C_2$. Then we consider a weighted geodesic ray transform
\beqq\label{eq-xray}
\msx^w f(z, z') = \int_0^{r(z, z')} \mcw(\gamma_{z, z'}(s))  f(\gamma_{z, z'}(s)) ds
\eeqq 
We note that $W$ is uniquely determined by the ray transform $\msx^wW$ and the stability estimate holds
\beqq\label{eq-stabest}
\|W\|_{L^2(\msm)} \leq C  \|\msx^w W\|_{H^1(\msc)}
\eeqq
where $C$ is a constant independent of $W$. This follows from  an application of Theorem \ref{thm-xraystab} to $\msx(\mcw W)$.   \\
 
Now we analyze $\mcf_{lead}$. With $\mcg_{0}(h, z, z') = (4\pi h)^{-2}e^{-i\sigma /h r(z, z')}|g(z')|^{-1/4}$, we have 
\beq 
\begin{gathered}
\mcf_{lead} =  h^2 \mcg_{0}(h)W \mcg_{0}(h) \\
= (4\pi)^{-4} h^{-2} \int_M e^{-i\sigma r(z, z')/h} |g(z)|^{-1/4} W(z') e^{-i\sigma r(z', z'')/h} |g(z')|^{-1/4}dz'    \\
  = (4\pi)^{-4} h^{-2} |g(z)|^{-1/4} \int_M e^{-i \frac{\sigma}{h}(r(z, z') + r(z', z''))}|g(z')|^{-1/4} W(z') dz'  
\end{gathered}
\eeq
At this point, we will apply stationary phase argument for a non-homogeneous phase function 
\beq
\Phi(z, z', z'') = r(z, z') + r(z', z''), \quad   z' \in \msm
\eeq
and we consider $z, z''\in \p \msm.$ 
To find critical points, we see from 
\beq
\p_{z'}\Phi(z, z', z'') = \p_{z'} r(z, z') + \p_{z'}r(z', z'') = 0
\eeq
 that $\p_{z'} r(z, z') = - \p_{z'}r(z', z'')$. This happens if and only if $z'$ is on $\gamma_{z, z''}$ the unique distance minimizing geodesic between $z, z'$. To see the rank of the  Hessian $\p_{z'}^2\Phi$, it is helpful to look at the Carleson-Sj\"olin condition in the estimates of oscillatory integrals, see \cite{MSS, So}. 
 
 Consider an oscillatory integral of the form
 \beq
 S_\la f(x) = \int_{\mbr^n} e^{i\la \phi(x, y)} a(x, y) f(y) dy
 \eeq
 where $a\in C_0^\infty(\mbr^n\times \mbr^n)$. For our problem, $\phi(x, y) = r(x, y)$ which is smooth away from $x = y$. The real valued smooth function $\phi$ satisfies the Carleson-Sj\"olin condition in this case if $\nabla_x \phi, \nabla_y \phi$ never vanish and 
 \beqq\label{eq-rank0}
 \text{rank} \phi_{xy}'' = n-1
 \eeqq
 and that there is a neighborhood $\mcu$ of $\supp a$ so that the immersed hypersurfaces 
 \beq
 \Sigma_{x_0} = \{\phi'_x(x_0, y): (x_0, y)\in \mcu\} 
 \eeq
 have everywhere non-vanishing Gaussian curvature. If $\phi = r(x, y)$, then 
 \beq
 \Sigma_{x_0} = \{\xi \in \mbr^n : \sum_{j, k = 1}^n g^{jk}(x_0)\xi_j\xi_k = 1\}
 \eeq 
 The curvature condition implies that 
 \beqq\label{eq-rank1}
\text{rank}(\frac{\p^2}{\p y_j \p y_k} \langle \phi'_{x}, \theta \rangle ) = n-1
 \eeqq
 where $\pm \theta \in \mbs^{n-1}$ is the  directions for which $\nabla_{y'} \langle \phi_{x}', \theta \rangle = 0.$ In fact, $\theta$ is orthogonal to $\Sigma_{x_0}$ at $\xi.$ For the distance function $r$ on $(\msm, g)$, the Carleson-Sj\"olin condition holds.
 
Let $z, z''\in \p \msm$ and $z'\in K$ a compact set of $\msm$. We observe that   
 \beq
 \begin{gathered}
\p_{z'} \Phi(z, z', z'') = \p_{z'} r(z, z') + \p_{z'}  r(z', z'') 
 =   \p_{z'} r(z, z') - \p_{z'} r(z'', z') 
 \end{gathered}
 \eeq
We let $\tilde z, \tilde z''\in \msm$ be points on a neighborhood of $z'$ such that 
\beq
 \p_{z'} r(z, z')  = \p_{z'} r(\tilde z, z'), \quad  \p_{z'} r(z'', z') = \p_{z'}r(\tilde z'', z')
\eeq
These are unit tangent vectors. Actually, if we let $\gamma_{z, z'}(s)$ be the geodesic with $\gamma_{z, z'}(0) = z', \dot \gamma_{z, z'}(0) = \p_{z'}r(z, z')$, then we can take $\tilde z = \gamma_{z, z'}(s_0)$ for some $s_0$ small. We can find $\tilde z''$ similarly on $\gamma_{z, z'}(s)$. 
Now we find that 
 \beq
 \begin{gathered}
\p_{z'}^2 \Phi(z, z', z'')  
 =  \p_{z'}( \p_{z'} r(\tilde z, z') - \p_{z'} r(\tilde z'', z') ) \approx \frac{\p^2}{\p z' \p z'}\p_z r(\tilde z, z')(\tilde z - \tilde z'')
 \end{gathered}
 \eeq
 for $\tilde z, \tilde z''$ close to $z'$ and $|\tilde z - \tilde z''|$ small. 
 Then the Carleson-Sj\"olin condition tells that the Hessian has rank $n-1 = 2$.  
Now one can choose local coordinates $z = (x, y), x\in \mbr, y = (y_1, y_2)\in \mbr^2$ such that the rank in \eqref{eq-rank1} in $y$ variable is $2$.  
We can perform stationary phase in $y$ variable and obtain using $h/|\sigma|$ as the small parameter. 

We recall the standard stationary phase expansion, see \cite[Proposition 1.2.4]{Dui}. Let $Q$ be a non-negative and symmetric matrix on $\mbr^2$ depending continuously on parameter $a\in \mbr^m$. Then
\beq
\begin{gathered}
\int e^{it \langle Q(a)y, y \rangle/2}g(y, a, t)dy \sim (\frac{2\pi}{t}) |\det Q(a)|^{-1/2} e^{\pi i \text{sgn}Q(a)/4} \\
\cdot \sum_{k = 0}^\infty \frac{1}{k!} (i \langle Q(a)^{-1} \p_y, \p_y\rangle /2 )^k g(0, a, t) t^{-k}
\end{gathered}
\eeq
for $t\rightarrow \infty$ uniformly in $a$. Applying this result, we get 
\beq 
\begin{gathered}
 \mcf_{lead}(z, z'') =  C h^{-2} |g(z)|^{-1/4} e^{i\sigma r(z, z'')/h} (\frac h\sigma)  \int_{0}^{r(z, z'')} |\det \text{Hess}\Phi(\gamma_{z, z''}(s))|^{-1/2} \\
 \cdot  |g(\gamma_{z, z''}(s))|^{-1/4}W(\gamma_{z, z''}(s)) ds + \mcf_{lead, 1} 
\end{gathered}
\eeq
where $C$ is a non-vanishing constant independent of $g, W$. For the remainder term, we have  
\beqq\label{lm-est0}
\begin{gathered}
|\mcf_{lead, 1}|_{C^0(\p \msm\times \p \msm)} \leq C (\frac h\sigma)^2 h^{-2}|W|_{C^2(\msm)},\\
|\mcf_{lead, 1}|_{C^1(\p \msm\times \p \msm)} \leq C (\frac h\sigma) h^{-2}|W|_{C^2(\msm)}.
\end{gathered}
\eeqq
Finally, we need the $\p_\nu$ derivative and we find that
\beqq\label{eq-fleadnu}
\begin{gathered}
 \p_\nu \mcf_{lead} = h^2 \p_\nu \mcg_{0}(h) W \mcg_{0}(h) \\
 =  C  h^{-2} \p_\nu r(z, z'') e^{i\sigma r(z, z'')/h} \msx^w(W) +  C  h^{-2} \frac{h}{\sigma} e^{i\sigma r(z, z'')/h} \p_\nu \msx^w(W)  + \p_\nu \mcf_{lead, 1} 
\end{gathered}
\eeqq
where $\p_\nu \mcf_{lead, 1}$ satisfy \eqref{lm-est0} as well.

\section{Estimates of the remainder terms}\label{sec-est}
First, we consider $\mcf_{para}$ and $\p_\nu \mcf_{para}$. 
\begin{lemma}\label{lm-est1}
Let $\sigma \in \mbc\backslash 0$ with $\im z\leq 0$. For $h < h_0$ small depending on $g$ and $|V|_{C^2(\msm)}$, the Schwartz kernel for $z, z'\in \p \msm$ satisfy 
\begin{enumerate}
\item $|r(z, z')\mcf_{para}(h, z, z')|_{C^0(\p \msm\times \p \msm)} \leq C |W|_{C^0(\msm)}$
\item $|r^2(z, z')\p_\nu \mcf_{para}(h, z, z')|_{C^0(\p \msm\times \p \msm)} \leq C h^{-1}| W|_{C^0(\msm)}$
\end{enumerate}
\end{lemma}
\bpf
We recall that 
\beq
 \mcf_{para} = -\mcg_0 h^2 W (\mcg_1 + \mcf) - (\mcg_1 + \mcf) h^2 W \mcg_0 - (\mcg_1 + \mcf) h^2 W (\mcg_1 + \mcf)
 \eeq
 For  term $\mcg_0 W \mcg_1$ and $\mcg_1 W\mcg_1$, we can apply the stationary phase argument as in Section \ref{sec-geo} to get the conclusion. For $\mcg_0 h^2 W \mcf$ and $\mcg_1 h^2 W \mcf$, one can estimate the integral directly as $\mcf = O(h^{-1})$ instead of $O(h^{-2})$.
 \epf

Next, we consider $\mcf_{res}, \p_\nu \mcf_{res}$. 
\begin{lemma}\label{lm-est2}
Let $\sigma \in \mbc\backslash 0$ with $\im z\leq 0$. For  $ h < h_0$ small depending on $g, |V|_{C^2(\msm)}$, the Schwartz kernel for $z, z'\in \p \msm$ satisfy 
\begin{enumerate}
\item $|r(z, z')\mcf_{res}(h, z, z')|_{C^0(\p \msm\times \p \msm)} \leq  C |W|^2_{C^0(\msm)}$
\item $|r^2(z, z')\p_\nu \mcf_{res}(h, z, z')|_{C^0(\p \msm\times \p \msm)} \leq C h^{-1} |W|^2_{C^0(\msm)} $
\end{enumerate}
\end{lemma}
\bpf
We start with the formula
\beq
\mcf_{res} = \mcr(h) h^2 W \mcr(h) h^2 W \mcr(h)(\id + h^2 W\mcr(h))^{-1}.
\eeq
We will study the kernel of $\mcf_{res, 1} = \mcr(h) h^2 W \mcr(h) h^2 W \mcr(h)$ 
because we know that the Schwartz kernel of $(\id + h^2 W\mcr(h))^{-1}$ is integrable and the $L^1$ norm is bounded in $h$ as $h\rightarrow 0.$  
We use that $\mcr = \mcg + \mcf$ and see that the Schwartz kernel is of the form 
\beq
\begin{gathered}
\mcf_{res, 1}(z, z''') = C h^{-6}\int_{M\times M} e^{-i\frac{\sigma}{h}r(z, z')} h^2 W(z') e^{-i\frac{\sigma}{h}r(z', z'')} h^2 W(z'')e^{-i\frac{\sigma}{h}r(z'', z'''')}\\
\cdot \frac{A(z, z', z'', z''')}{r(z', z'')}dz'dz''
\end{gathered}
\eeq
where $z \in \p \msm$ and $z', z'', z'''\in \msm$ and the amplitude $A$ is smooth.  
We would like to apply the stationary phase argument but the distance function is not smooth at the diagonal. Let's consider the phase function $\Phi(z, z', z'') = r(z, z') + r(z', z'')$ with integration in $z'\in \msm$. Here, $z\in \p \msm, z''\in \msm$. Since $W$ is supported away from $\p \msm$, we just need to consider when $z'$ is close to $z''$. 
For fixed $z''\in \msm$, we let $B_\eps(z'')$ be the ball of radius $\eps$ centered at $z''$ and we split the integral (for fixed $z'',z'''$)
\beq
\begin{gathered}
 \int_\msm e^{-i\frac{\sigma}{h}r(z, z')} h^2 W(z') e^{-i\frac{\sigma}{h}r(z', z'')} \cdot \frac{A(z, z', z'', z''')}{r(z', z'')}dz' \\
 =  \int_{\msm\backslash B_\eps(z'')} e^{-i\frac{\sigma}{h}r(z, z')} h^2 W(z') e^{-i\frac{\sigma}{h}r(z', z'')}  
\cdot \frac{A(z, z', z'', z''')}{r(z', z'')}dz' \\
+ \int_{B_\eps(z'')} e^{-i\frac{\sigma}{h}r(z, z')} h^2 W(z') e^{-i\frac{\sigma}{h}r(z', z'')}  
\cdot \frac{A(z, z', z'', z''')}{r(z', z'')}dz'  = I_1 + I_2
\end{gathered}
\eeq
For integral $I_1$, we can apply stationary phase argument as before to conclude that $|I_1(z, z'')|\leq Ch^3|W|_{C^0}$. For $I_2$, we change the integral to polar coordinate $z' = z'' + \rho w, \rho\in (0, \eps), w\in \mbs^{1}$ and get 
\beq
\begin{gathered}
I_2 = \int_{\mbs^1}\int_0^\eps e^{-i\frac{\sigma}{h}r(z, z'' + \rho w)} h^2 W(z'' + \rho w) e^{-i\frac{\sigma}{h}\rho}  
\cdot \frac{A(z, z''+\rho w, z'', z''')}{\rho}\rho^2 d\rho dw\\
 =  \int_{\mbs^1}\int_0^1 e^{-i\frac{\sigma}{h}r(z, z'' + s h w)} h^2 W(z'' + sh w) e^{-i \sigma s}  
 A(z, z'' +sh w, z'', z''') sh hds dw
\end{gathered}
\eeq
where we changed variables $s = \rho/h$ in the second line. Thus $|I_2(z, z'')| \leq C h^4|W|_{C^0}$.  
 To summarize, we get 
 \beq
\begin{gathered}
\mcf_{res, 1}(z, z''') = C h^{-3}\int_\msm e^{-i\frac{\sigma}{h}r(z, z'')} \mca_0(z, z'')  h^2 W(z'')e^{-i\frac{\sigma}{h}r(z'', z'''')} 
\cdot A(z, z', z'', z''')  dz''
\end{gathered}
\eeq
 where $|\mca_0|_{C^0} \leq C |W|_{C^0}$ and is smooth in $z, z''$. Now we can apply the stationary phase argument again to get  
\beq
|\mcf_{res, 1}(z, z''')|_{C^0(\p \msm\times \p \msm)} \leq  C |W|^2_{C^0(\msm)}
\eeq
Similarly, we get 
\beq
|\p_\nu \mcf_{res, 1}(z, z''')|_{C^0(\p \msm\times \p \msm)} \leq  Ch^{-1}  |W|^2_{C^0(\msm)}
\eeq
This completes the proof of the lemma. 
\epf

\section{The semiclassical DtN map}\label{sec-dtn}
In this section, we look for the dependency of the DtN map on $h$. 
For the classical Calder\'on problem, it is known, see for example  \cite{LU}, that  the DtN map is a pseudo-differential operator of order one on $\p \msm$. The approach there is to  decompose the elliptic operators in boundary normal coordinates. The method implicitly relies on standard elliptic regularity results which can be studied in this approach, see   \cite{Tre1}. For the semiclassical problem, we re-examine the approach and pay attention to the dependency on $h.$  We carry out the construction for dimensions $n\geq 2$ as the argument is the same.

We recall the decomposition of $\lap_g$ from  \cite[Section 2]{LU}. Consider the Laplace-Beltrami operator $\lap_g$ in a boundary normal coordinates $(x_1, \cdots, x_{n})$ in which $\p \msm = \{x_n = 0\}$  and we denote $x' = (x_1, \cdots, x_{n-1}).$ The metric is of the form 
\beq
g = \sum_{\alpha, \beta = 1}^{n-1}g_{\alpha\beta} dx^\alpha dx^\beta + (dx^n)^2
\eeq
We recall that 
\beq
\lap_g = - \sum_{i, j = 1^n}  \delta^{-\ha} \p_{x_i} (\delta^\ha g^{ij} \p_{x_j} )
\eeq
where $\delta = |\det g_{ij}|$. We write $D_{x_n} = -i \p_{x_n}$. In boundary normal coordinates, we have  
\beq
\lap_g = D_{x_n}^2 + i E(x) D_{x_n} + Q(x, D_{x'})
\eeq
where
\beq
\begin{gathered}
E(x) = -\ha \sum_{\alpha, \beta} g^{\alpha\beta}(x)\p_{x_n}g_{\alpha\beta}(x)\\
Q(x, D_{x'}) = \sum_{\alpha, \beta = 1}^n g^{\alpha\beta} D_{x_\alpha}D_{x_\beta} - i \sum_{\alpha, \beta} (\ha g^{\alpha\beta}(x) \p_{x^\alpha} \log \delta(x) + \p_{x^\alpha}g^{\alpha \beta})) D_{x^\beta}
\end{gathered}
\eeq
see \cite[page 1101]{LU}. 
It is proved in \cite[Proposition 1.1]{LU} that there exists a pseudodifferential operator $A(x, D_{x'})$ of order $1$ in $x'$ depending smoothly on $x_n$ such that 
\beq
\lap_g = (D_{x_n} +i E(x) - i A(x, D_{x'}))(D_{x_n} + i A(x, D_{x'})) + B
\eeq
where $B$ denotes a smoothing operator. Now we consider the problem
\beq
(\lap_g + V)u - \la^2 u = 0
\eeq
We use the decomposition to get
\beq
(D_{x_n}  +i E(x) - i A(x, D_{x'}))(D_{x_n}  + i A(x, D_{x'}))u  + B u - \la^2 u   = 0
\eeq
Dividing by $\la^2$ and setting $h = 1/\la$, we get 
\beq
(hD_{x_n}  +i hE(x) - i hA(x, D_{x'}))(hD_{x_n}  + i hA(x, D_{x'}))u  + (h^2B - 1) u   = 0
\eeq
Now we convert the elliptic problem  to a system 
\beqq\label{eq-heateq}
\begin{gathered}
(hD_{x_n}   + i hA(x, D_{x'}))u = v\\
(hD_{x_n}  +ihE(x) - i hA(x, D_{x'}))v  + (h^2B - 1) u   = 0
\end{gathered}
\eeqq
The boundary condition $u = f$ at $\p \msm$ is converted to $u = f$ at $x_n = 0.$  Notice that $h A(x, D_{x'})$ is a semiclassical pseudo-differential operator of order $1$. Later, we will denote it by $A(x, h D_{x'})$ to signify the semiclassical nature. 

Let's see what we need to do to find the DtN map. We let $U(x_n, h)$ be a semiclassical parametrix such that 
\beq
\begin{gathered}
(hD_{x_n}   + i A(x, hD_{x'}))U(x_n, h) = 0 \text{ mod } O(h^\infty) \\
U(0, h) = \id 
\end{gathered}
\eeq
Then we can write using Duhamel's principle 
\beqq\label{eq-para1}
u(x_n, x') = U(x_n, h) f(x') + \int_0^{x_n} U(x_n - x_n', h)v(x_n, x') dx_n'
\eeqq
modulo a $O(h^\infty)$ term. To find $v$, 
we let $W(x_n, h)$ be the parametrix of the backward heat equation 
\beq
\begin{gathered}
(hD_{x_n}   + ih E -  i A(x, hD_{x'}))W(x_n, h) = 0 \text{ mod } O(h^\infty)  \\
W(T, h) = \id, \quad T > 0.
\end{gathered}
\eeq
We can write 
\beqq\label{eq-para2}
v(x_n, x') = W(x_n) v(T, x') - \int_{x_n}^T W(x_n - x_n')(h^2B  + 1) u(x_n', x')dx_n'
\eeqq
modulo a $O(h^\infty)$ term. 
The important property of the parametrix we need to establish is that for $t > 0$, they are smoothing operators of order $h^\infty.$ We know a priori that the solution $u$ to the Dirichlet problem is smooth and the $H^1$ norm is bounded by $h^{-2}$ which can be seen from the variational form. Using \eqref{eq-para2} and choosing smooth $v(T, x')$ bounded in $h^{-2}$, we conclude that $v$ is smooth and is of order $h^\infty.$ The choice of $v(T, x')$ will not change the regularity of $u$, just as in page 131 of \cite{Tre1}.  Finally, using \eqref{eq-para1}, we obtain solution $u(t, x)$ up to a $h^\infty$ error term which is also smooth. The $h^\infty$ is important because a classical construction does not improve the order of $h.$ 

After this is done, we obtain the DtN map from the first equation of \eqref{eq-heateq}
\beq
D_{x_n} u =   h^{-1} i A(x, hD_{x'})u   = h^{-1} A(x, hD_{x'})f \text{ at } x_n = 0
\eeq
up to a $h^\infty$ smooth term.  We will prove that
\begin{prop}\label{prop-dtn}
Let $(M, g)$ be as in Theorem \ref{thm-main} except that the dimension of $M$ is $n\geq2$. 
Consider the semiclassical Dirichlet problem 
\beq 
\begin{gathered}
h^2(\lap_g + V) - u  = 0  \text{ in } \msm, \\
u = f \text{ at } \p \msm
\end{gathered}
\eeq  
 For the DtN map $\La$ defined in \eqref{eq-dtn1}, the Schwartz kernel $\La(z, z'), z, z'\in \p M$ is such that $|z - z'|{\La}(z, z')$ is continuous and bounded for $h$ small. Here, $|\cdot|$ denotes Euclidean norm. 
\end{prop}
The key in this approach is to construct a parametrix of the semiclassical heat equation which we study next. 

\section{Parametrix of semiclassical heat equations}\label{sec-heat}
We brief review the basics of semiclassical quantization from \cite{Zw}. 
For $h\in [0, 1)$, consider $a(h, x, \xi) \in C^\infty([0, 1); S^m(\mbr^n_x; \mbr^n_\xi))$. Here, $S^m(\mbr^n; \mbr^n)$ is the standard symbol class, which is the set of $C^\infty$ functions on $\mbr^n_x \times \mbr^n_\xi$ satisfying 
\beq
|D_x^\alpha D_\xi^\beta a(h, x, \xi)| \leq C_{\alpha \beta} \langle \xi\rangle^{m - |\beta|}
\eeq
for all $\alpha, \beta \in \mbn^n$. The estimate is uniform on compact set of $\mbr_x^n$. When the context is clear, we also abbreviate the notation as $S^m.$ The semiclassical operator with symbol $a$ is defined as 
\beq
A(x, hD) u(x) =  \op_\hbar(a)u(x) = (2\pi h)^{-n} \int_{\mbr^n}e^{i(x - y)\cdot \xi/h} a(h, x, \xi) u(y)dy d\xi
\eeq
Here, we only use the standard quantization. The semiclassical principal symbol is $\sigma_{\hbar, m}(A) = a|_{h = 0} \in S^m$.

For $h\in [0, 1), t\geq 0, x', \xi'\in \mbr^{n-1}$, we consider symbols $a(h, t, x', \xi') \in C^\infty([0, 1)_h\times [0, \infty); S^m(\mbr^{n-1}_{x'}, \mbr^{n-1}_{\xi'})$. We assume that 
\beq
a(h, t, x', \xi') \sim \sum_{j = 0}^\infty h^j a_j(t, x', \xi'), \quad a_j \in C^\infty([0, 1)_h; S^m). 
\eeq
Moreover, the semiclassical principal symbol $a_0$ is elliptic, namely 
\beq
\text{$|a_0(t, x', \xi')|\geq \gamma |\xi'|^m$ for $\gamma > 0$ and all $(t, x')$. }
\eeq
Let $A(t, x', hD)$ be the semiclassical quantization  of $a(h, t, x', \xi')$.  We consider semiclassical heat equations  
\beqq\label{eq-heat}
\begin{gathered}
(h D_t + A(t, x', hD_{x'}) )u = 0, \text{ for } t > 0, x'\in \mbr^{n-1}\\
 u = f, \quad t = 0
\end{gathered}
\eeqq
The argument below follows closely Section 1.1 of Chapter III of \cite{Tre1}. 
We aim to find a semiclassical parametrix of the form 
\beq
U(t, x', hD_{x'})  = \sum_{j = 0}^\infty h^j U_j(t, x', hD_{x'})
\eeq
where $U_j$ are semiclassical quantizations of $u_j(t, x', \xi') \in S^{-m}$, that is 
\beq
U_j(t)f(x') = (2\pi h)^{-(n-1)} \int e^{i(x' - y')\xi'/h} u_j(t, x', \xi')  f(y') d\xi' dy'
\eeq
We write $x = (t, x') \in \mbr^n$ below. Formally, we have
\beqq\label{eq-U1}
\begin{gathered}
(h D_t + A(x, hD_{x'}) )U(x, h D_{x'}) \\
=  \sum_{j = 0}^\infty h^{j } h \p_t U_j(x, hD_{x'}) +  \sum_{j = 0}^\infty h^{j} A(x, hD_{x'}) U_j(x, hD_{x'}) \\
\end{gathered}
\eeqq
Let $C_j(x, hD_{x'}) = A(x, hD_{x'})U_j(x, hD_{x'})$ be the composition. We use semiclassical calculus for standard quantization to conclude that $C_j$ are semiclassical pseudo-differential operators with full symbols 
\beq
 c_j(h, x, \xi') = \sum_{k = 0}^N \frac{h^k}{k!} (i\langle D_{\xi'}, D_{y'}\rangle)^k (a(x', \xi') u_j(y', \eta'))|_{y' = x', \eta' = \xi'} + O(h^{N+1})
\eeq
as $h\rightarrow 0.$ Here, $c_j\in C^\infty([0, 1)_h; S^0).$ See Theorem 4.14 and Theorem 4.18 of \cite{Zw}. 
Using the asymptotic expansion of $a(x, \xi)$, we find that 
 the symbol expansion of $A(x, h D_{x'})U(x, hD_{x'})$ is  $\sum_{j = 0}^\infty h^j d_{j}(t, x', \xi'),$ where 
\beq
\begin{gathered}
 d_{j}(t, x', \xi') 
= \sum_{k = 0}^j \sum_{l = 0}^k \frac{1}{l!} (i\langle D_{\xi'}, D_{y'}\rangle)^{l} (a_{k - l}(x', \xi') u_{j - k}(y', \eta'))|_{y' = x', \eta' = \xi'}.  
\end{gathered}
\eeq
Using these formula, we get from \eqref{eq-U1} that
\beqq\label{eq-U2}
\begin{gathered}
(h D_t + A(x, hD_{x'}) )U(x, hD_{x'})=  (2\pi h)^{-(n-1)} \cdot \\
 \sum_{j = 0}^\infty h^{j } \int_{\mbr^{n-1}}  e^{i(x' - y')\xi'/h}  
\cdot (h \p_t u_j(t, x', \xi')  +  d_{j}(t, x', \xi')) d\xi' 
\end{gathered}  
\eeqq
From the  order $h^j$ terms, we get equations
\beqq\label{eq-symu1}
h \p_t u_j(t, x', \xi') +  a_0(t, x', \xi')u_j(t, x', \xi')  + e_j(t, x', \xi') = 0, \quad \forall j\geq 0
\eeqq
where $e_0 = 0$ and for $j\geq 1$ we have
\beq
\begin{gathered}
e_j(t, x', \xi') = d_j(t, x', \xi') -  a_0(t, x', \xi')u_j(t, x', \xi') \\
= \sum_{k = 1}^j \sum_{l = 0}^k \frac{1}{l!} (i\langle D_{\xi'}, D_{y'}\rangle)^{l} (a_{k - l}(x', \xi') u_{j - k}(y', \eta'))|_{y' = x', \eta' = \xi'}
\end{gathered}
\eeq
Notice that the term $e_j$ involves $u_k, k < j$. So we will solve these equations iteratively. 
The  equation \eqref{eq-symu1} comes with initial conditions. At $t = 0$, we get 
\beq
U_0(0, x', hD_{x'}) = \id, \quad U_j(t, x', hD_{x'}) = 0,  \quad j\geq 1. 
\eeq
This implies that 
\beqq\label{eq-symu2}
u_0(0, x', \xi') = 1, \quad u_j(0, x', \xi') = 0, \quad j\geq 1. 
\eeqq

To solve  \eqref{eq-symu1} with \eqref{eq-symu2}, we look for solutions of the form
\beq
u_j(t, x', \xi') = (2\pi i)^{-1} \int_{\gamma} e^{\rho z t/h} k_j(t, x', \xi', z) dz 
\eeq
where $\rho  = \langle \xi' \rangle^{m}$, $z\in \mbc$ is a complex parameter and $\gamma$ is a contour so that the integrand is holomorphic in a neighborhood of $\gamma$ for $t, x'$ in a compact set.   
Plug this into \eqref{eq-symu1}, we obtain
\beq
\begin{gathered}
(2\pi i)^{-1} \int_\gamma e^{\rho z t/h} \mcl k_j dz = 0,\\
\mcl k_j =  \rho z k_j + h D_t k_j + a_0 k_j + e_j
\end{gathered}
\eeq
Now $hD_t k_j$ has an additional $h$. So we should look at the asymptotics in \eqref{eq-U2}. After re-arrangment, we get from the order $h^j$ term the equations 
\beq
\begin{gathered}
\tilde \mcl k_0 =  \rho z k_0 + a_0 k_0 \\ 
\tilde \mcl k_j =  \rho z k_j + a_0 k_j + e_j + D_t k_{j-1}, \quad j\geq 1
\end{gathered}
\eeq
We aim to solve $k_j$ from $\tilde \mcl k_j = \rho$. Then we justify that the choices solve equations \eqref{eq-symu1}, \eqref{eq-symu2}. 

First, we recall that $a_0 \in S^m$ is elliptic. Thus $\rho^{-1}a_0 \in S^0$ is also elliptic and we see that for $\xi' \in \mbr^{n-1}$, $z - \rho^{-1} a_0$ 
is non-zero for $\im z < 0. $ Let $E(z) = (z - \rho^{-1}a_0)^{-1}$.  
For $j = 0$,  we solve $\tilde\mcl k_0 = \rho$ to get $k_0 = E(z)$.  
For $j\geq 1$, we get 
\beq
k_j = - E(z)\rho^{-1} [D_t k_{j-1} + e_j] 
\eeq
It follows from the argument in Page 137 of \cite{Tre1} that $k_j$ are smooth in $(t, z)$ and valued in $S^{m_j}$ and $m_j \leq -j \min(1, m).$ 
This finishes the construction.  Finally, we prove the regularizing properties of the parametrix. 
\begin{lemma}
The $u_j$ defined above are smoothing and belong to $O(h^\infty)$. 
\end{lemma}
\bpf
We look at 
\beq
u_j(t, x', \xi') = (2\pi i)^{-1}\int_{\gamma} e^{\rho z t/h} k_j(t, x', \xi', z) dz 
\eeq
in which $k_j \in C^{\infty}([0, \infty); S^{m_j})$. We  estimate
\beq
|\p_{x'}^\alpha \p_{\xi'}^\beta \p_t^\gamma k_j(t, x', \xi')| \leq C \langle \xi' \rangle^{m_j - |\beta|}. 
\eeq
Also, we estimate
\beq
\begin{gathered}
|\p_t^\gamma \p_{\xi'}^\beta e^{\rho z t/h}| \leq C (t/h)^{-N}(1 + |\xi'|)^{-N}
\end{gathered}
\eeq
Putting the estimate together, we obtain that $u_j$ is indeed smoothing and of $h^\infty.$ 
\epf

Now we can finish the proof of Proposition \ref{prop-dtn}. 
\bpf[Proof of Proposition \ref{prop-dtn}] 
 For the DtN map $\La$, we know that modulo a smooth $O(h^\infty)$ term, the Schwartz kernel is  
\beq
\begin{gathered}
\La(x', y')   =  h^{-1} (2\pi i h)^{-(n-1)} \int e^{i(x' - y')\cdot \xi'/h} a(x', \xi', h) d\xi'\\
 = h^{-1}(2\pi i)^{-{n-1}}\int e^{i(x' - y')\cdot \xi'} a(x', h\xi', h) d\xi'
\end{gathered}
\eeq
Because $a$ is a symbol of order $1$, the Schwartz kernel has a singularity like $|x' - y'|$  and $|x' - y'|\La(x', y')$ is bounded in $h.$
\epf

\section{Proof of Theorem \ref{thm-main}}\label{sec-proof}
Let's consider two potentials $V, \tilde V$ on $\msm$ and $\La, \tilde \La$ be the corresponding DtN maps.  Suppose $\La = \tilde \La.$ We  use \eqref{eq-ide} 
\beq
\begin{gathered}
 \p_\nu \tilde \mcr(z, z', h) - \p_\nu \mcr(z, z', h) 
 = \int_{\p \msm} \tilde \mcr(z, z'', h) \tilde\La(z'', z') dz''   - \int_{\p \msm} \mcr(z, z'', h)  \La(z'', z') dz''  
\end{gathered}
\eeq
and  \eqref{eq-approx}
 \beq
 \begin{gathered}
 \mcr(z, z', h) - \tilde \mcr(z, z', h)  
     =  \mcf_{lead} + \mcf_{para} + \mcf_{res} 
     \end{gathered}
 \eeq
and the estimates  of the remainder terms to finish the proof. 

We start from the left hand side of \eqref{eq-ide}
\beqq\label{eq-final1}
\begin{gathered}
 \p_\nu \tilde \mcr(z, z', h) - \p_\nu \mcr(z, z', h)   
 =  C  h^{-2} \p_\nu r(z, z'') e^{i\sigma r(z, z'')/h} \msx^w(W) \\
 +  C  h^{-2} \frac{h}{\sigma} e^{i\sigma r(z, z'')/h} \p_\nu \msx^w(W)  + \p_\nu \mcf_{lead, 1}  +\p_\nu \mcf_{res} + \p_\nu \mcf_{para}.  
\end{gathered}
\eeqq
For the right hand side of \eqref{eq-ide}, we have
\beqq\label{eq-ide3}
\begin{gathered}
 \int_{\p \msm} \tilde \mcr(z, z'', h) \tilde\La(z'', z') dz''   - \int_{\p \msm} \mcr(z, z'', h)  \La(z'', z') dz'' \\ =  \int_{\p \msm} (\tilde \mcr(z, z'', h) -  \mcr(z, z'', h) ) \tilde \La(z'', z') dz'' + 
 \int_{\p \msm}   \mcr(z, z'', h) (\tilde \La(z'', z') - \La(z'', z')) dz''
\end{gathered}
\eeqq
Since $\La = \tilde \La$, we only need to consider the first term. But we remark that the second term will give the stability estimate, although we do not pursue it here. We have 
\beqq\label{eq-final2}
\begin{gathered}
 \int_{\p \msm} (\tilde \mcr(z, z'', \sigma, h) -  \mcr(z, z'', \sigma, h) ) \tilde \La(z'', z') dz''  \\
 = \int_{\p \msm} (C h h^{-2} |g(z)|^{-1/4} e^{i\sigma r(z, z'')/h}\msx^w W \tilde \La(z', z'') dz' dz''  \\
 + \int_{\p \msm}(\mcf_{lead, 1}(z, z')  + \mcf_{res}(z, z') + \mcf_{para}(z, z'))\tilde \La(z', z'')dz''. 
\end{gathered}
\eeqq
Now we can use the estimate of $\tilde \mcr - \mcr$ and the kernel estimate of $\tilde \La$ in Proposition \ref{prop-dtn}. 
Also, we will use the stability estimate of $\msx.$ First of all, 
\beq
\begin{gathered}
\|\msx^wW\|_{L^2(\msm)} \leq Ch^2 (h^{-1}\|\p_\nu \msx^wW\|_{L^2(\msc)} + | \p_\nu \mcf_{lead, 1}  +\p_\nu \mcf_{res} + \p_\nu \mcf_{para}|_{C^0(\p \msm\times \p \msm)}) \\
+ C h |\msx^w W|_{C^0(\p\msm\times\p\msm)} + Ch^2  |\mcf_{lead, 1}  + \mcf_{res} + \mcf_{para}|_{C^0(\p \msm\times \p \msm)}
\end{gathered}
\eeq
The estimates of these terms are done in Lemma \ref{lm-est1}, \ref{lm-est2} and \eqref{lm-est0}. Also, we used that the boundary $\p \msm$ is strictly convex and that we only consider $z, z'\in \p \msm$ and away from $z = z'$ in \eqref{eq-ide}, so the absolute value of the $\p_\nu r$ term in \eqref{eq-final1} is bounded from below.  
For $\|\p_\nu\msx^w W\|_{L^2(\msc)}$, we just need to notice that 
\beq
\begin{gathered}
\p_{z^i} \msx^w W = \p_{z^i} \int_0^{r(z, z')} \mcw(\gamma_{z, z'}(s))W(\gamma_{z, z'}(s)) ds 
 = \int_0^{r(z, z')} \p_{z^i}\mcw(\gamma_{z, z'}(s))W(\gamma_{z, z'}(s)) ds
\end{gathered}
\eeq
because $W$ is compactly supported. 
Thus 
\beq
\|\msx^w W\|_{H^1(\msc)} \leq  \|\msx^w W\|_{L^2(\msc)} + \|\msx^w (\p W)\|_{L^2(\msc)}.
\eeq
Using these estimates, we obtain from \eqref{eq-final1}, \eqref{eq-final2} that 
\beq
h^{-2}\|W\|_{L^2(\msm)} \leq  C h^{-2} h |W|_{C^2(\msm)} + Ch^{-1}|W|^2_{C^0(\msm)}
\eeq
so 
\beqq\label{eq-ine}
\|W\|_{L^2(\msm)} \leq  Ch |W|_{C^2(\msm)} + C h |W|^2_{C^0(\msm)} \leq Ch |W|_{C^2(\msm)}  
\eeqq
if  $|W|_{C^0(\msm)}\leq C_0$. 

To prove Theorem \ref{thm-main}, we use $|W|_{C^2(\msm)}\leq C_1 \|W\|_{L^2(\msm)}$ for some $C_1$. 
Then we conclude that $\|W\|_{L^2(\msm)} \leq Ch \|W\|_{L^2(\msm)}$. So for $h$ sufficiently small, $W = 0$ which  completes the proof.

 \begin{remark}\label{rmk-stab}
For a class of potential, it is possible to obtain some stability estimate.  Let $K$ be a compact set of $\msm$. For $C_0, C_1 >0$, we define 
\beqq\label{def-vset}
\begin{gathered}
\msv = \{V\in C_0^\infty(\msm):  \supp V \subset K, |V|_{C^2(\msm)}\leq C_0,  
 |V|_{C^2(\msm)}\leq C_1 \|V\|_{L^2(\msm)}\}.
 \end{gathered}
\eeqq
For potentials in $\mcv$, it follows from the proof of Theorem \ref{thm-main} above that such potentials are uniquely determined by their DtN map for a sufficiently small $h$ which depends on $C_0, C_1$. One can obtain H\"older type stability estimate by examining the last term in \eqref{eq-ide3}.  This agrees with he phenomena of increased stability for high frequency Schr\"odinger operators on $\mbr^n$, which has been studied in the literatures.  See the recent work \cite{KUW} and the references therein.

Also, we remark that the last inequality in \eqref{def-vset} resembles the so-called inverse inequality in numerical methods, see for example \cite[Section 6.2]{LaTh}. 
\end{remark}

 \begin{remark}\label{rmk-recon}
 Our proof leads to an approximate reconstruction method. From \eqref{eq-final1}, we get that 
\beq
\msx^w (\tilde V - V) = \frac{Ch^2}{\p_\nu r(z, z')} \int_{\mbr} \mcr(z, z'', h) (\tilde \La(z'', z') - \La(z'', z')) dz'' + O(h)
\eeq
where the weight and the constant can be found explicitly from the proof. Note that they only depend on the background manifold $(\msm, g)$. Thus, to reconstruct a potential $\tilde V$ from $\tilde \La$, we can choose a reference potential $V = 0$ with corresponding $\La$ which can be computed for  the  manifold $(\msm, g)$. 
Therefore, for $h$ sufficiently small, we can find $\tilde V$ by inverting the weighted geodesic ray transform.
\end{remark}

\section*{Acknowledgement}
GU was partly supported by NSF, a Walker Family Endowed Professorship at University of Washington and a Si-Yuan Professorship at IAS, HKUST.


\end{document}